   \documentclass[centertags]{amsart}
   \newtheorem{lemma}{Lemma}[section]
   \newtheorem{theorem}[lemma]{Theorem}
   \newtheorem{remark}[lemma]{Remark}
   
   \newtheorem{coro}[lemma]{Corollary}
   \newtheorem{definition}[lemma]{Definition}

 \newcommand{\p}{\partial}

   \title
   [The Stochastic 3D Quasigeostrophic Flows]
   {The 3D  Quasigeostrophic  Fluid Dynamics \\under Random
   Forcing on Boundary}

 \author{Jinqiao Duan}
   \address[Jinqiao Duan]
   {Corresponding author: Department of Applied Mathematics\\
   Illinois Institute of Technology\\
   Chicago, IL 60616, USA.\\
Tel: (312)567-5335 $\;\;$ Fax: (312)567-3135}
   \email
   [Jinqiao Duan]{duan@iit.edu}

   \author{Bj{\"o}rn Schmalfu{\ss} }
   \address[Bj{\"o}rn Schmalfu{\ss}]
   {Department of Applied Sciences\\
   University of Technology and Applied Sciences\\
   Geusaer Stra{ss}e\\
   D--06217 Merseburg\\
   Germany\\}
   \email
   [Bj{\"o}rn Schmalfu{\ss}]{schmalfuss@in.fh-merseburg.de}

   \date{December 15,  2000 }

   \subjclass{Primary 60H25, 76D05, 47H10; Secondary 34D35, 86A05}
   \keywords{\, 3D Quasigeostrophic flows,  wind forcing on
 boundary,  stable stationary flows, random dynamical systems,
random attractors}

\begin{document}

\begin{abstract}
The three-dimensional baroclinic quasigeostrophic flow model has
been widely used to study basic mechanisms in oceanic flows and
climate dynamics. In this paper, we consider this flow model under
random wind forcing and time-periodic fluctuations on fluid
boundary (the interface between the oceans and the atmosphere).
The time-periodic fluctuations are due to periodic rotation
of the earth and thus periodic exposure of the earth to the solar
radiation.
After establishing the well-posedness of the baroclinic
quasigeostrophic flow model in the state space, we demonstrate the
existence of the random attractors, again in the state space.  We
also discuss the relevance of our result
 to climate modeling.

\bigskip
{\bf PACS Codes}: 02.50.Fz, 47.32.-y, 92.10.Fj, 92.60.Gn

\end{abstract}

 \maketitle


\newpage

   \section{Introduction}

We will study the quasigeostrophic equation (QGE)
for the ocean in a three dimensional domain.  This geophysical
flow model has been formally
 derived as an approximation of
the rotating three-dimensional primitive   equations \cite{Ped87}.
Bourgeois and Beale  \cite{BeaBou94}, Embid and Majda
\cite{EmbMaj96}, and Desjardins  and  Grenier \cite{DesGre98} have
recently shown that the three-dimensional   quasigeostrophic
equation is a valid approximation of the   primitive equations in
the limit of zero Rossby number. Holm \cite{Hol86}  has
established the Hamiltonian formulation for the inviscid
quasigeostrophic equation. In particular, we will consider the
following version \cite{Wan92, Dua99} of the equation in terms of
the stream function $\psi(x,y,z, t)$:
\begin{equation*}
\begin{split}
\tilde\Delta
   \psi_t&+J(\psi,\tilde\Delta\psi)+\beta\psi_x=\nu\tilde\Delta \tilde\Delta\psi\\
& \;\;\;
 \mbox{with} \;\;\;
\tilde\Delta \psi
=\psi_{xx}+\psi_{yy}+\left(\frac{f_0^2}{N^2(z)}\psi_z\right)_z
\end{split}
\end{equation*}
Here  $x,\,y,\,z$ are  Cartesian coordinates in zonal (east),
meridional (north), vertical directions, respectively;
 $f_0 +\beta y$ (with $f_0, \beta$ constants) is the Coriolis
parameter;  $N(z)>0$ is the Brunt-Vaisala stratification
frequency, and $ \nu > 0$ is viscosity.  Moreover, $J(f, g)
=f_xg_y -f_yg_x$ is the Jacobi operator and  potential vorticity
is defined as $ \tilde\Delta \psi+ f_0+ \beta y$. Note that
$\tilde\Delta\psi$  can be regarded  as  a modified Laplacian
operator  where the coefficient in the vertical $z$ direction is
adjusted due to the density stratification, and the coefficients
in $x, y$ directions are constants due to the  horizontal density
homogenity in the 3D quasigeostrophic flow model formulation.
Bennett and Kloeden \cite{BenKlo81}  have also used a similar
modified Laplacian viscous term in a more complicated 3D
quasigeostrophic flow model involving thermodynamics as well as
 hydrodynamics.
\\
Our aim is to study the potential verticity evolution in an ocean
under the influence of the atmosphere.  Following an idea of
Hasselmann  \cite{Has76}  one can divide the geophysical  or
meteorological flow  into two parts. These two parts are the
slowly changing climate part  and rapidly changing weather part.
The weather part can be modeled by a stochastic process such as
white noise, see Hasselmann \cite{Has76},  Arnold \cite{Arn00} and
Saltzman \cite{Sal83}. Oceanic flows are affected (on the ocean
surface) by these short time influences due to weather variations
which are usually called wind  forcing. Moreover, oceanic flows
are also affected  by climatic variations due to periodic rotation
of the earth and thus periodic exposure of the earth to the solar
radiation; see \cite{PeiOor92}, Chapter 6 and \cite{Ler98},
Chapter 11.

Since the exchange between the atmosphere and an ocean takes place
at the surface of the ocean,  we will consider the above
quasigeostrophic partial differential equation with white noise
Neumann boundary condition on the top surface of the ocean
\cite{Ped87}.  We
will also add some time-periodic boundary condition on the ocean
surface due to  periodic  rotation of the earth and  the solar
radiation.    Since there is no influence of the weather at the
bottom of the ocean we will have the homogeneous boundary
condition there. The boundary conditions in horizontal directions
are assumed to be periodic  for mathematical  convenience as in
other recent works \cite{EmbMaj96, BeaBou94}.  However,  since we
have a forth order differential operator at the right hand side of
the quasigeostrophic equation, we need a second group of boundary
conditions, say,  for $\tilde \Delta \psi$, and these conditions
are the ones used in \cite{Wan92}, \cite{Dua99} to be specified
below. For simplicity we will suppose that this second group
 of boundary
conditions are deterministic.  But generalizations are possible.\\

Our aim is to find structures in the dynamics of the QGE.  As we
have time dependent random and time-periodic boundary conditions,
we will obtain a nonautonomous dynamical system with random
influences. We will
 show how to find attractors for such a dynamical system.
The existence and interpretation of climatic attractors have been
controversial and have caused a lot of debate \cite{Nicolis}. A
low dimensional climatic attractor was regarded as an indication
that the main feature of long-time climatic evolution may be
viewed as the manifestation of a {\em deterministic} dynamics.
Our result is about {\em random} attractors,  and thus the long
time regimes that such attractors may represent still carry the
stochastic information of the geophysical flow system.
Stochastically forced QGE has been used to investigate various
mechanisims in geophysical flows \cite{Holloway, Mul96,
Samelson, BraJinWan98}.

The QGE can be transformed into an evolution equation with
standard boundary conditions. For this transformation we need an
Ornstein-Uhlenbeck process fulfilling our dynamical random or
time-periodic boundary conditions. This transformation will be
introduced  in Section 2.\\
In  Section 3 we investigate the coefficients of the
 transformed
evolution equation, and further obtain a  global
existence and uniqueness result and some   regularity
result. \\
In Section 4, we study the random dynamics
of the transformed QGE.
Based on   the uniqueness result above,  the  transformed evolution
equation generates a nonautonomous dynamical system.  In addition,
if we restrict this system  to discrete time step of the period of
the periodic  rotation of the earth, we obtain a random dynamical
system. This random dynamical system has a random attractor. This
result can be
extended to the dynamical system on the real time axis.\\
The Section 5 contains the proofs.
 Finally, we summarize our results
in Section 6.


   \section{The 3D quasigeostrophic equation}

   Let $O=(0,2\pi)^3$ be the  cube which is a model for a piece of the
   ocean. For $x,\,y,\,z\in (0,2\pi)$
   and smooth functions $u(x,y,z),\,v(x,y,z)$ we define the Jacobi operator
   \begin{equation*}
   J(u,v)=u_xv_y-u_yv_x.
   \end{equation*}
   In addition, the differential operator $\tilde\Delta$
   is defined by
   \begin{equation*}
  \tilde\Delta := \p_{xx}+\p_{yy}+\left(F(z)\p_z\right)_z
   \end{equation*}
   for $f_0\in\mathbb{R}\setminus{0},\,F(z):=\frac{f_0^2}{N^2(z)}$
and $N(z)$ is defined to be a positive
   $C^\infty$-smooth function $N(z)>0$ on $[0,2\pi]$.
   Let $\nu,\,\beta$
   be positive constants.
   In the following we investigate the 3D QGE flow model
   \cite{Ped87, Sal98, Wan92}:
   \begin{equation}\label{eq-0.5}
   \tilde\Delta
   \psi_t+J(\psi,\tilde\Delta\psi)+\beta\psi_x=\nu\tilde\Delta\tilde\Delta\psi.
   \end{equation}
   We impose the following boundary conditions for this equation.
   Let $O_{\cdot,\cdot,0},\,O_{\cdot,\cdot,2\pi},\cdots$ be the faces of the
   cube
   $O$ then
   we assume periodic boundary conditions in $x,\,y$ directions
   \[
   \psi|_{O_{0,\cdot,\cdot}}=\psi|_{O_{2\pi,\cdot,\cdot}}\qquad
   \psi_x|_{O_{0,\cdot,\cdot}}=\psi_x|_{O_{2\pi,\cdot,\cdot}}
   \]
   and similar for the  faces $O_{\cdot,0,\cdot},\,O_{\cdot,2\pi,\cdot}$
together with derivative
   in $y$ direction.
   With respect to the $z$ direction, $\psi$   fulfills
   Neumann boundary conditions on top and bottom boundaries;
see \cite{Ped87}, p.367.  Namely,  we impose that
   \[
   \frac{\partial \psi}{\partial n}=\psi_z=0
 \quad \text{on  bottom boundary} \; O_{\cdot,\cdot,0}
   \]
   and
   \[
   \frac{\partial \psi}{\partial n}=\psi_z=\dot{W}\quad \text{on top boundary} \;
  O_{\cdot,\cdot,2\pi} \;\;,
   \]
   where $\dot{W}$ is a temporal white noise and $n$ denotes the unit outer normal vector.\\
   Moreover, $\tilde\Delta\psi$ is supposed to be periodic in $x, y$
   directions:
   \[
   \tilde\Delta\psi|_{O_{0,\cdot,\cdot}}=\tilde\Delta\psi|_{O_{2\pi,\cdot,\cdot}}\qquad
   \tilde\Delta\psi_x|_{O_{0,\cdot,\cdot}}=\tilde\Delta\psi_x|_{O_{2\pi,\cdot,\cdot}}
   \]
   and similar for the  faces $O_{\cdot,0,\cdot},\,O_{\cdot,2\pi,\cdot}$ with derivative
   in $y$ direction.
   For the other faces  on the top and the bottom of the ocean,
we impose homogeneous Neumann
   boundary conditions as in \cite{Wan92, Dua99, BenKlo81}
   \[
   \tilde\Delta\psi_z|_{O_{\cdot,\cdot,0}}=0,\qquad
   \tilde\Delta\psi_z|_{O_{\cdot,\cdot,2\pi}}=0
   \]

   In addition,  we impose that
   \[
   \int_O\psi dO=0,
   \]
and
\[
   \int_O\tilde\Delta\psi dO=0.
   \]
We also assume an appropriate initial condition
\[
\psi(x, y,z, 0) = \psi_0(x,y,z).
\]

Later on we will see that we can find
an Ornstein-Uhlenbeck process
$\eta$ satisfying
the linear differential equation
\begin{equation}\label{eq0}
    \eta_t= \nu \tilde \Delta \eta ,
   \end{equation}
where $\eta$ fulfills the same boundary conditions
as $\psi$ in some sense.

   We now transform formally the above equation to another parabolic
   differential equation
   with random coefficients but with standard boundary conditions.
   To do this we have to introduce the solution operator $G$ of the following
   elliptic problem:
   \begin{equation}\label{eq1}
   \tilde\Delta\xi=f
   \end{equation}
   on $O$ with periodic boundary conditions in $x, y$ directions and
   {\em homogeneous} Neumann boundary conditions in $z$ direction.
   For existence and regularity properties of this operator see below.
   $G$ can be considered as a linear continuous operator. Here $f$
   is an element in a particular Hilbert space.
   \\
   For the difference of $\psi$ and $\eta$ we obtain for
   (\ref{eq-0.5})
   \begin{eqnarray*}
   (\tilde\Delta\psi-\tilde\Delta\eta)_t&+&J(\psi-\eta,\tilde\Delta\psi-\tilde\Delta\eta)+
   J(\eta,\tilde\Delta\psi-\tilde\Delta\eta)+J(\psi-\eta,\tilde\Delta\eta)
   +J(\eta,\tilde\Delta\eta)\\
   &+&\beta(\psi-\eta)_x+\beta\eta_x=\nu\tilde\Delta\tilde\Delta (\psi-\eta)
   \end{eqnarray*}
   where $\eta$ is assumed to be sufficiently regular.
   The difference $\psi-\eta$ satisfies the same boundary condition as
   $\xi$
   in (\ref{eq1}). Hence by the formal substitution
   $v=\tilde\Delta\psi-\tilde\Delta\eta$
   we can write
   \begin{align*}
   v_t&+J(G(v),v)+J(\eta,v)+J(G(v),\tilde\Delta\eta)+J(\eta,\tilde\Delta\eta)\\
   &+\beta G(v)_x+\beta\eta_x=\nu\tilde\Delta v.
   \end{align*}
   The problem (\ref{eq0}) satisfies only standard boundary conditions
   with respect  to $\eta$ but not with respect to $\tilde\Delta\eta$.
   To obtain an equation with standard boundary conditions
   we need the transformation $u=v+\tilde\Delta\eta$ assuming formally that $\eta$ is
   sufficiently regular.
   Then $u$ fulfills periodic boundary conditions in $x, y$ directions,
and satisfies the following random evolution:
   \begin{equation}\label{eq2}
   u_t+J(G(u),u)+J(\eta-G(\tilde\Delta\eta),u)+\beta
   G(u)_x-\beta(G(\tilde\Delta\eta))_x-\eta_x)=\nu\tilde\Delta u.
   \end{equation}
   This equation will be treated in the rest of the paper.\\
   \begin{remark}
   {\rm
   To obtain a solution of the original equation all reverse
   transformations make sense if $\eta$ is contained in a Sobolev
   space of second order. Hence we can find a solution of this
   equation contained in this Sobolev space.  Thus, in the following
we will formulate our results for this transformed random
evolution equation (\ref{eq2}).\\
Note that $\eta-G(\tilde\Delta\eta)$ can be handled as a function
fulfilling the boundary conditions with
$\tilde\Delta(\eta-G(\tilde\Delta\eta))=0$.
   }
   \end{remark}


   \section{Existence and uniqueness of the solution for the
                   3D quasigeostrophic equation }

In this section we consider the well-posedness of the
transformed QGE (\ref{eq2}).
   We intend to consider equation (\ref{eq2}) as an evolution equation
   on a rigged space $V\subset H\subset V^\prime$.
   Since this equation is given with respect to periodic boundary
   conditions in horizontal directions the space $H$ is defined to be
   \[
   \dot{L}_2(O)=\{u\in L_2(O),\;\int_OudO=0\}.
   \]
   The inner product on $H$ is denoted by $(\cdot,\cdot)$.
   Let $H^k,\,k\in\mathbb{N}$ be the usual Sobolev  space consisting of functions
   with square integrable derivative up to $k$-th order.
   If $k$ is not an integer these spaces are the usual Slobodeckij
   spaces.
   For $V\subset H^1\cap H$ we choose the space of  functions from $H^1$ such
   that
   \[
   u|_{O_{0,\cdot,\cdot}}=u|_{O_{2\pi,\cdot,\cdot}},
   \qquad
   u|_{O_{\cdot,0,\cdot}}=u|_{O_{\cdot,2\pi,\cdot}}.
   \]
   Note that for functions in $H^1$ the trace on the boundary is
   well defined. This set will be equipped with the usual $H^1$ inner product
   denoted by $(\cdot,\cdot)_V$.
   Let $A$ be a linear bounded operator
   \[
   A:V\to V^\prime
   \]
   which is the usual operator stemming from the bilinear form $a(u,v)$
   defined by $-\tilde \Delta$ with periodic boundary conditions in $x, y$
   directions and homogeneous Neumann boundary condition in $z$ direction.
   Note that $-\tilde \Delta$ is symmetric with respect to these
   boundary conditions.
   We have that
   \[
   \langle Au,u\rangle=\|u\|_V^2.
   \]
   $\langle\cdot,\cdot\rangle$ denotes the dual pairing between
   $V^\prime=D(A^{-\frac{1}{2}})$ and $V=D(A^{\frac{1}{2}})$.
   The operator $A$ is an isomorphism from $V$ to $V^\prime$.
   Since $V$ is compactly embedded in $H$ we can
   define $G(f)=A^{-1}f$ which is a continuous operator for instance from $H$
   to $D(A)$ or more generally  from $D(A^s)$ to $D(A^{s+1}),\,s\in\mathbb{R}$.
   For the definition of the spaces $D(A^s)$ and their norms
   $\|\cdot\|_{D(A^s)}$ see Temam \cite{Tem97} Section II.2.
    In particular,
   for $f\in H$ the function $G(f)$ is periodic in $x,y$
   directions.\\
   On $H$ the operator $A$ has the spectrum
   $0<\lambda_1\le\lambda_2\le \cdots$. The associated
   eigenfunctions are denoted by $e_1,\,e_2,\cdots$.\\
   The embedding constant between $V$ and $H$ is given by $\lambda_1$:
   \[
   \lambda_1\|u\|_H^2\le \|u\|_V^2.
   \]

   We now investigate the properties of the operator $J$ which defines the
   nonlinearity of (\ref{eq2}).
   \begin{lemma}\label{l1}
   Suppose that $m_1,\,m_2,\,m_3$ are three nonnegative numbers
   less than $\frac{3}{2}$ with
   \[
   \sum_{i=1}^3m_i\ge \frac{3}{2}.
   \]
   Then there exists a constant $c_1>0$ such that for $u\in
   H^{m_1+1},\,v\in H^{m_2+1}$ and
   $w\in H^{m_3}$
   \[
   |\langle J(u,v),w\rangle|\le c_1\|u\|_{H^{m_1+1}}
   \|v\|_{H^{m_2+1}}\|w\|_{H^{m_3}}.
   \]
   \end{lemma}
   \begin{proof}
   On account  the H{\"o}lder
   inequality, we get
   \[
   |\langle J(u,v),w\rangle|\le \|\nabla u\|_{L_{q_1}}\|\nabla
   v\|_{L_{q_2}}\|u\|_{L_{q_3}},
   \qquad
   \frac{1}{q_1}+\frac{1}{q_2}+\frac{1}{q_3}\le 1.
   \]
   The embedding property  $H^{m}\subset
   L_q(O)$ where $\frac{1}{q}=\frac{1}{2}-\frac{m}{3}$
   for $q>1,\,m\ge 0,\,m\not=\frac{3}{2}$ gives the conclusion. For the idea of the
   proof see Temam \cite{Tem83},  Section 2.3.
   \end{proof}
   \begin{remark}\label{xx1}
   {\rm
   Suppose    two of the $m_i$'s in the last lemma, say
   $m_2,\,m_3$, have the value zero. Then  if $m_1$  is chosen
   bigger than $\frac{3}{2}$ the conclusion of the last lemma
   remains true. This follows if we apply at first the Sobolev lemma
   and then the Cauchy-Schwarz inequality.
   }
   \end{remark}
   From Lemma \ref{l1}  we can derive some  a priori estimates
   for the nonlinearity of (\ref{eq2}).
   \begin{coro}\label{c0}
   Suppose that $u\in D(A^\alpha),\,\alpha>\frac{3}{4}$ and $v,\,w\in V$. Then
   \[
   |\langle J(u,v),w\rangle|\le
   c_1\|u\|_{D(A^\alpha)}\|v\|_V\|w\|_V,
   \]
   and hence for $u\in H,\,v,\,w\in V$
   \[
   |\langle J(G(u),v),w\rangle|\le c_1\|u\|_{H}\|v\|_V\|w\|_V.
   \]
   Suppose now that $u\in D(A^\alpha),\,\alpha>\frac{5}{4}$ and $v\in V,\,w\in H$.
   Then
   \[
   |\langle J(u,v),w\rangle|\le
   c_1\|u\|_{D(A^\alpha)}\|v\|_V\|w\|_H,
   \]
   and hence for  $u\in D(A^\alpha),\,\alpha>\frac{1}{4}$ and $v\in V,\,w\in H$
   \[
   |\langle J(G(u),v),w\rangle|\le
   c_1\|u\|_{D(A^\alpha)}\|v\|_V\|w\|_H.
   \]
   Furthermore, if $u\in H^2$ and $v,\,w\in V$ we have the inequality
   \[
   |\langle J(u,v),w\rangle|\le c_1\|u\|_{H^2}\|v\|_V\|w\|_V.
   \]
   \end{coro}
   Based on this Corollary, we also get
   \begin{lemma}\label{l4}
   For any $\mu>0$ there exists a constant $c_{2,\mu}$ such that for
   $v\in H$ and $u\in H^2$
   \[
   |\langle J(u,v),Av)\rangle|\le c_{2,\mu}\|u\|_{H^2}^2\|v\|_V^2+\mu\|Av\|_H^2.
   \]
   and
   for an $\eta\in H^2$
   \[
   |\langle J(\eta,v),Av)\rangle|\le
   c_{2,\mu}^\prime\|\eta\|_{H^2}\|v\|_{D(A^\frac{3}{4})}\|Av\|_0\le
   c_{2,\mu}\|\eta\|_{H^2}^4\|v\|_V^2+\mu \|Av\|_H^2
   \]
   \end{lemma}
   For the last inequality we need an interpolation argument and the Young inequality.
   Later on we can use this lemma to derive additional regularity
   properties of (\ref{eq2}).

   Also on account of Lemma \ref{l1} we obtain the following
   algebraic properties of $J$.

   \begin{lemma}\label{c1}
   If $v,\,w\in V$ and $u\in H^2$ and periodic in $x$, $y$ directions,
   then we have
   \[
   \langle J(u,v),w\rangle=-\langle J(u,w),v\rangle.
   \]
   Hence we have $\langle J(u,v),v\rangle=0$.
   \end{lemma}
   \begin{proof}
   We choose $u$ from a set of sufficiently smooth functions which is dense in $D(A^\alpha)$.
   Integration by parts formulae yields
   \begin{align*}
   \begin{split}
   \int_0^{2\pi}\int_{O_{\cdot,\cdot,z}}&u_xv_ywdxdydz-\int_0^{2\pi}\int_{O_{\cdot,\cdot,z}}u_yv_xwdxdydz\\
   &=
   -\int_0^{2\pi}\int_{O_{\cdot,\cdot,z}}u_{xy}vwdxdydz+\int_0^{2\pi}\int_{O_{\cdot,\cdot,z}}u_{yx}vwdxdydz
   \\
   &
   -\int_0^{2\pi}\int_{O_{\cdot,\cdot,z}}u_xvw_y
   dxdydz+\int_0^{2\pi}\int_{O_{\cdot,\cdot,z}}u_yvw_xdxdydz
   \\
   &
   +\int_0^{2\pi}\int_0^{2\pi}u_xvw\big|_{y=0}^{y=2\pi}dxdz
   -\int_0^{2\pi}\int_0^{2\pi}u_yvw\big|_{x=0}^{x=2\pi}dydz.
   \end{split}
   \end{align*}
   Note that the two boundary terms are zero. Indeed, since
   $u(x,0,z)=u(x,2\pi,z)$, we know that $u_x(x,0,z)=u_x(x,2\pi,z)$, and since $v,\,w$
   are
   $2\pi$--periodic with
   respect to $x$. By the smoothness assumption we can suppose that the derivatives
   on the boundary are well defined.
   Thus the first
   boundary term is zero.
   Similarly, the second boundary term is also zero.
   By the continuity of $J$ (see Corollary \ref{c0}) we can extend
   the relation to $H^2\times V\times V$.
   The second
   claim is due to antisymmetric property of the first claim of this corollary.
   \end{proof}

   The following lemma will be used to obtain the continuity of
   the solution operator of (\ref{eq2}).
   \begin{lemma}\label{l3}
   There exists a constant $c_3>0$ such that for $u_1,\,u_2\in V$
   \[
   |\langle J(G(u_1),u_1)-J(G(u_2),u_2),u_1-u_2\rangle|\le
   c_3\|u_1-u_2\|_V\|u_1-u_2\|_H\|u_1\|_V.
   \]
   \end{lemma}
   \begin{proof}
   By Lemma \ref{c1} the left hand side above expression is equal to
   \begin{align*}
   |\langle J(G(u_1),u_1),u_2\rangle&+\langle J(G(u_2),u_2),u_1\rangle|\\
   &
   =|\langle J(G(u_1),u_1),u_2-u_1\rangle-\langle
   J(G(u_2),u_1),u_2-u_1\rangle|\\
   &=|\langle J(G(u_1)-G(u_2),u_1),u_1-u_2\rangle|.
   \end{align*}
   Corollary \ref{c0} gives the conclusion.
   \end{proof}

  The equation  (\ref{eq2}) can be considered as an evolution
equation on the rigged space
   $V\subset H\subset V^\prime$ introduced at the beginning of this section.\\
   We now explain the properties of coefficients which are contained in (\ref{eq2}).
   Due to Corollary
   \ref{c0} we have a bilinear continuous operator
   \[
   B(\cdot,\cdot):V\times V\to H,\quad B(u,v):= J(G(u),v).
   \]
   In addition, we have a time dependent linear operator
   \[
   C(t,u):= J(\eta(t)-G(\tilde\Delta \eta(t)),u),\quad C(t):V\to V^\prime
   \]
   where $\eta(t)\in L_{2,loc}(0,\infty;H^2\cap H)$ and periodic in $x,y$ directions.
   Then $G(\tilde\Delta \eta(t))\in D(A)$.

   \begin{lemma}\label{ll}
   Suppose that $\eta(\cdot)\in L_{2,loc}(0,\infty;H^2)$
    and  $\eta$ is periodic in $x$,$y$
   directions.
   Then we have for almost any $t\ge 0$:
   \[
   \|C(t,\cdot)\|_{V,V^\prime}\le c_4\|\eta(t)\|_{H^2},\qquad
   \langle C(t,u),u\rangle=0.
   \]
   \end{lemma}
   Indeed,  because $\eta(t)-G(\tilde\Delta\eta(t))\in H^2$ for any $t\ge 0$ we obtain
   the first conclusion by Corollary \ref{c0}. The second Conclusion follows
from Lemma \ref{c1}.\\

   We now investigate the last linear operator appearing in (\ref{eq2})
   which is defined by
   \[
   D(u)= \beta G(u)_x:V\to H.
   \]
   We obtain straightforwardly
   $\|D(\cdot)\|_{V,H}\le c_5$.
   \begin{lemma}\label{ll1}
   For $u\in V$ we have
   \[
   (D(u),u)=0
   \]
   \end{lemma}
   \begin{proof}
   Denoting $G(u)\in D(A)$ by $\xi$ we have
   \begin{align*}
   \beta^{-1}(D(u),u)=&(\xi_x,\tilde\Delta
   \xi)=\frac{1}{2}\bigg(\int_O(\xi_x^2)_xdxdydz
   -(\xi_y^2)_xdxdydz\\
   -&\int_0^{2\pi}\int_0^{2\pi}
   F(z)\int_0^{2\pi}(\xi_z^2)_xdxdydz\bigg)
   \end{align*}
   via the integration by parts. The second term under the integral generates the
   boundary term
   \[
   \frac{1}{2}\int_0^{2\pi}\int_0^{2\pi}\xi_x\xi_y\big|_{y=0}^{y=2\pi}dxdz.
   \]
   Indeed, for sufficiently smooth $\xi$ from a dense set in $V$ we have
   $\xi(x,0,z)=\xi(x,2\pi,z)$  and thus  $\xi_x(x,0,z)=\xi_x(x,2\pi,z)$. It follows
   from the periodicity
   in $y$ we have $\xi_y(x,0,z)=\xi_y(x,2\pi,z)$ such that this boundary term is zero.\\
   For the last term the following
   boundary
   term appears
   \[
   \frac{1}{2}\int_0^{2\pi}\int_0^{2\pi}
   F(z)\xi_x\xi_z\bigg|_{z=0}^{z=2\pi}dxdy
   \]
   which is zero by the homogeneous Neumann boundary conditions.
   Integration with respect to $x$ and using the periodicity
   as in the proof of Lemma \ref{c1}, we get the conclusion.
   \end{proof}
   We can formulate the evolution equation  (\ref{eq2}) on $V\subset H\subset V^\prime$
   \begin{equation}\label{eq2.2}
   u_t+\nu Au+B(u,u)+C(t,u)+D(u)=f(t),\quad u(0)=u_0\in H
   \end{equation}
   for $f(t)=\beta(G(\tilde\Delta \eta(t))_x-\eta(t)_x)$ which is contained in
   $L_{2,loc}(0,\infty;H)$.

   Apart from the linear operators $C(t)$ and $D$ equation (\ref{eq2.2})
   has the form  of  equations of 2D Navier -Stokes type; see Temam
   \cite{Tem79},
   Chapter 3,  for which we have existence and uniqueness.
Here {\em 2D} means that the
   conclusion of Lemma \ref{l3} is fulfilled which is responsible for a
   uniqueness theorem. By the regularity
   properties
   of the  operators $B,\,C(t)$ and $D$ the same method as for the 2D
   Navier-Stokes equations ensures existence
and uniqueness for (\ref{eq2.2}). Thus we get the following
main result in this section, about the well-posedness
for 3D quasigeostrophic flows under random
wind forcing on ocean surface.

   \begin{theorem}\label{tex}
(Well-posedness)
   Suppose that the Ornstein-Uhlenbeck process $\eta$ defined in (\ref{eq0})
is in  $L_{2,loc}(0,\infty;H^2)$.  Then the 3D QGE
equation (\ref{eq2.2})
 or   (\ref{eq2})
   has a unique solution $u(t)\in L_{2,loc}(0,\infty;V)\cap C([0,\infty);H)$.
   \end{theorem}
   \begin{remark}\label{rneu}
   {\rm
   i) On account of Lemma \ref{l3} we can prove that $u(t)$
   depends continuously on the initial conditions $x\in H$.
   This will be used later on.\\
   ii) Due to Lemma \ref{l4} for $t>0$ and a bounded set of
   initial conditions in $H$ the image with respect to the solution operator
   $u_0 \to u(t)$ is bounded in $V$. Hence the solution  operator  is compact for
   $t>0$.
   }
   \end{remark}


   \section{The dynamical system of 3D quasigeostrophic flows}

In this section we study the dynamical behavior of QGE
(\ref{eq2.2}).
   In the following we are going to describe the background perturbations
   defined on the boundary of $O$ which will influence the dynamical system
   generated by (\ref{eq2.2}).
   We will have two different influences. The first perturbation is a
   white noise which models the weather or the small scale impact of
the atmospheric motion through wind forcing
 on the surface.   The other one is a
   periodic motion which serves as a model for the  impact due to
 periodic rotation of the
   earth and thus periodic exposure of the earth
to the solar radiation; see \cite{PeiOor92}, Chapter 6 and \cite{Ler98},
Chapter 11.
\\
   In the first part of this section we are going to explain a dynamical model of the
   boundary conditions.
   \\
   We consider the elliptic differential equation
   \begin{equation*}
   \tilde\Delta u=0,\quad\frac{\partial u}{\partial n}=0\;
   \text{on } O_{\cdot,\cdot,0}\quad
   \frac{\partial u}{\partial n}=f\in
   H^{k+\frac{1}{2}}(O_{\cdot,\cdot,2\pi})\cap
   \dot{L}_2(O_{\cdot,\cdot,2\pi}).
   \end{equation*}
   On the other faces we have periodic boundary conditions.
   This equation has a unique solution; see Egorov and Shubin \cite{EgoShu91} Page 130f.
   The solution operator
   is a linear continuous operator with the image in $H^{k+2}$. We denote this operator
    by $\tilde G(f)$. $H^{k+\frac{1}{2}}(O_{\cdot,\cdot,2\pi})$ denotes a usual boundary
    space.
   \\

   At first we consider the random part of the boundary conditions. Let $W$ be a
   continuous temporal Wiener process with values in a
   Hilbert space
   $U=H^{2+\frac{1}{2}}(O_{\cdot,\cdot,2\pi})\cap \dot{L}_2(O_{\cdot,\cdot,2\pi})$
   This Wiener process is defined for positive and negative times; see Arnold
   \cite{Arn98} Page 547.
   The covariance operator of $W$ is denoted by $Q$, which is a
   positive symmetric linear operator on $U$ with finite trace. The dynamics
   of $W$ is given by the {\em metric dynamical system} consisting
   of a probability space $(\Omega,\mathcal{F},\,\mathbb{P})$ and a
   flow $\theta$,
   $(\Omega,\mathcal{F},\mathbb{P},\theta)$ where $\mathbb{P}$ is the
   Wiener measure with covariance $Q$ and
   $\theta=(\theta_t)_{t\in\mathbb{R}}$ is the flow of the Wiener shift:
   \[
   W(\cdot,\theta_t\omega)=W(\cdot+t,\omega)-W(t,\omega)\quad\text{for }
   t\in\mathbb{R}.
   \]
   The mapping $\theta$ is
   $(\mathcal{B}(\mathbb{R})\otimes\mathcal{F},\mathcal{F})$--measurable
   and fulfills the property
   \begin{equation}\label{eqo6}
   \theta_{t+\tau}=\theta_t\circ\theta_\tau,\qquad t,\,\tau\in\mathbb{R}.
   \end{equation}
   For instance, we can choose  $\Omega$ to be the set of continuous functions
   $C_0(\mathbb{R},U)$ which are zero at zero and for $\mathcal{F}$ we choose the
   Borel-$\sigma$-algebra of $C_0(\mathbb{R},U)$. Note that $\mathbb{P}$ is ergodic with respect to
   $\theta$. Later on we have to restrict this metric dynamical system to fix particular
   dynamical properties.\\
   We now show the existence of a solution of (\ref{eq0}) satisfying
   particular properties. At first we show that the following problem
   has a {\em stationary} solution:
   \begin{equation}\label{eqo1}
   {\eta_1}_t=\tilde\Delta \eta_1,\quad \frac{\partial\eta_1}{\partial
   n}=0,\quad (x,y)\in O_{\cdot,\cdot,0},
   \quad \frac{\partial\eta_1}{\partial n}=\dot{W}(t),\quad (x,y)\in
   O_{\cdot,\cdot,2\pi}
   \end{equation}
   under periodic boundary conditions with respect to the other faces of the
   cube
   $O$ and a initial condition at time $t=0$. This solution will serve as a process which
   compensates the nonhomogeneous boundary conditions in
   (\ref{eq2}).
   A similar problem has been considered in Da Prato and Zabczyk
   \cite{DaPZab96} Chapter 13 or \cite{DaPZab93}.
   \\
   In contrast to  (\ref{eqo1}) we now consider boundary conditions
   which are defined to be {\em homogeneous} Neumann boundary conditions on
   $O_{\cdot,\cdot,2\pi}$.
   The solution of (\ref{eqo1}) with these boundary conditions  generates an
   analytic $C_0$-semigroup  $(S(t))_{t\ge 0}$.
   This semigroup has a generator, which will be denoted
   also by $-\nu A$ where $A$ is equivalent to the operator introduced in
   the last section.\\
   The following lemma allows us to define a stationary
Ornstein-Uhlenbeck process which fulfills particular regularity
   assumptions.
   We notice that
$\tilde G$ maps continuously
$H^\frac{5}{2}(\mathcal{O}(\cdot,\cdot,2\pi))$ into $V\cap H^4$
where
\[
V=\{u \in H^1\cap H: u(0,\cdot,\cdot)=u(2\pi,\cdot,\cdot),\;
u(\cdot,0,\cdot)=u(\cdot,2\pi,\cdot)\}.
\]
We also note that for our boundary conditions $V\cap H^2\not=
D(A)$.
\begin{lemma}\label{lv1}
Let $Q$ be a linear operator of finite trace on  $U=H^\frac{5}{2}(O_{\cdot,\cdot,2\pi})\cap H$.
We consider the elliptic differential operators with respect to the boundary
conditions formulated in (\ref{eq1})
\begin{equation}\label{abslast}
\tilde\Delta u+(F_zu)_z,\qquad
\left\{
\begin{array}{l}
\tilde\Delta u+(F_zu)_z\\
\tilde\Delta v+2(F_zv)_z+F_{zz}v+F_{zzz}u
\end{array}
\right.
\end{equation}
which are supposed to be generators of analytic $C_0$-semigroups of negative type.
Then we have
\[
\int_0^\infty\|AS(\tau)\|_{\mathcal{L}_2(U,H^2)}^2<\infty,\qquad
\int_0^\infty\|AS(\tau)\|_{\mathcal{L}_2(U,D(A^\beta))}^2<\infty,\quad\beta<\frac{1}{4}.
\]
\end{lemma}
\begin{proof}
Let $\tilde e_i$ be a complete orthonormal base in $U$.
Then $\tilde GQ^\frac{1}{2}\tilde e_i\in H^4\cap H$.
We now study the $D_x,\,D_{xx},\cdots D_{zz}$-derivative
of $\tilde GQ^\frac{1}{2}\tilde e_i$. All these expressions are periodic
in $x,\,y$ directions. In addition, all these elements are
contained in $D(A^\eta),\,\eta<\frac{3}{4}$
This follows because we only have periodic
and Neumann boundary conditions,
 see Da Prato and Zabczyk
\cite{DaPZab92} Page 401.
Since $S(t)g$ is the solution of a differential
linear parabolic equation
\begin{equation}\label{eqev2}
\frac{dv}{dt}+\nu Av=0,\qquad v(0)=g.
\end{equation}
we can differentiate the
the solution with respect to $D_x$.
The coefficients of
the equation are independent of $x$. The function $D_xS(t)g$ fulfills the same
equation as $S(t)g$ but with derived initial condition. Hence we have  $D_xS(t)g=S(t)D_xg$.
The regularity of $D_xg\in D(A^\eta)$ ensures that
$\|AD_xS(t)g\|_H$ is well defined such that $D_xS(t)g\in D(A)$, see \cite{DaPZab92} Page 392:
\begin{equation}\label{ev1}
\|AD_xS(\tau)g\|_H^2=\|AS(\tau)D_xg\|_H^2\le \frac{C^2}{t^{2-2\eta}}e^{-2\alpha t}\|D_xg\|_{D(A^\eta)}^2
\quad t>0.
\end{equation}
for $g=\tilde GQ^\frac{1}{2}\tilde e_i$, $\alpha>0$ and $\eta\in (\frac{1}{2},\frac{3}{4})$.
But on $D(A)$ we
can identify $A$ with $-\tilde\Delta$ such that we can exchange $A$ and $D_x$.
Integration of (\ref{ev1}) from 0 to $\infty$ and summation over $i$ yields the
desired estimate.
Similarly, we can argue for the derivatives $D_x,\,D_y,\,D_{yy},\cdots$.\\
Similar considerations allow us to prove the estimates for
derivatives $D_{zx},\cdots,D_{zz}$. In contrast to the above
formulae a derivative $D_z$ changes $A$ to an operator $\hat A$
given in (\ref{eqev2}) with semigroup $\hat S(t)$. The assumption
(\ref{abslast} ensures that $\hat S(t)$ is of negative type with
our usual boundary conditions. We get
\begin{equation*}
\|D_zAS(\tau)g\|_H^2=
\|\hat AD_zS(\tau)g\|_H^2=\|\hat A\hat S(\tau)D_zg\|_H^2
\le \frac{\hat C^2}{t^{2-2\eta}}e^{-2\hat\alpha t}\|D_zg\|_{D(\hat A^\eta)}^2,\,t>0.
\end{equation*}
The same method gives us  estimates in $z$-direction. Note that if
we differentiate the coefficients of (\ref{eqev2}) twice in z
direction one term without a second $z$ derivative appears. To
take this term into account we need the system of differential
operators in (\ref{abslast}).

To obtain the second conclusion we use the estimate for $\beta\in (0,\frac{1}{4})$
\[
\|AS(t)g\|_{D(A^\beta)}^2=\|A^\beta AS(t)g\|_H^2
\le\frac{C^2}{t^{2-2\eta}}e^{-2\alpha t}\|g\|_{D(A^{\beta+\eta})}^2.
\]
and an appropriate $\eta>\frac{1}{2}$: $\eta+\beta\in (\frac{1}{2},\frac{3}{4}),\quad t>0$.
\end{proof}
Note that if $N(z)$ (hence $F(z)$) is constant then the assumption about the
generator is superflow.
\begin{remark}\label{remev3}
{\rm
On account of Da Prato and Zabczyk \cite{DaPZab96} Chapter 13
there exists an Ornstein-Uhlenbeck process
\[
\eta_1(\cdot,\omega)=A\int_0^\cdot S(\cdot-\tau)d\tilde G
W(\tau,\omega)\in L_{2,loc}(0,\infty;H^2)
\]
solving (\ref{eqo1}) with initial condition
zero at $t=0$. $\eta_1(t,\omega)$ is periodic in $x,\,y$ directions.
}
\end{remark}
Now we are going to show that $\eta_1(\cdot,\omega)$ has
continuous paths in $H$.
\begin{lemma}\label{lv2}
There exists a continuous version in
$H$ and in particular,
\[
\mathbb{E}\sup_{s\in[0,1]}\|\eta_1(s)\|_{H}^2<\infty.
\]
\end{lemma}
\begin{proof}
Due to Da Prato and Zabczyk \cite{DaPZab92} Theorem 5.9 and Remark
5.11 we have to check
\[
\int_0^\infty\tau^{-\gamma}\|AS(\tau)\tilde GQ^\frac{1}{2}
\|_{\mathcal{L}_2(U,H)}^2d\tau<\infty
\]
for some $\gamma>0$. The left hand side is equal to
\begin{equation}\label{eqv2}
\sum_{i=1}^\infty\lambda_i^{2}\|Q^\frac{1}{2}\tilde G^\ast
e_i\|_H^2 \int_0^te^{-2\lambda_i\tau}\tau^{-\gamma}d\tau,
\end{equation}
see Da Prato and Zabczyk \cite{DaPZab96} Proposition 13.2.4 where $(e_i)$ is orthonormal
base introduced above. For
sufficiently small $\gamma>0$ we have a $c(\gamma)>0$:
\[
\lambda_i^{-2\beta}\int_0^t\tau^{-\gamma}e^{-2\lambda_i\tau}d\tau
\le c(\gamma)\int_0^te^{-2\lambda_i\tau}d\tau.
\]
Hence (\ref{eqv2}) can be estimated by
\begin{align*}
\sum_{i=1}^\infty&\lambda_i^{2}\|Q^\frac{1}{2}\tilde
G^\ast e_i\|_H^2 \int_0^t\tau^{-\gamma}e^{-2\lambda_i\tau}d\tau\\
&\le c(\gamma)\sum_{i=1}^\infty\lambda_i^{2+2\beta}\|Q^\frac{1}{2}\tilde
G^\ast e_i\|_H^2 \int_0^te^{-2\lambda_i\tau}d\tau=
c(\gamma)\int_0^t\|AS(\tau)\|_{\mathcal{L}_2(U,D(A^\beta))}^2d\tau
\end{align*}
which is finite by Lemma \ref{lv1}.
\end{proof}
We now are looking for a stationary process solving (\ref{eqo1}).
\begin{lemma}\label{ll3}
There exists a random variable $\eta_1\in H$ such that
$t\to\eta_1(\theta_t\omega)$ solves (\ref{eqo1}).
\end{lemma}
\begin{proof}
By the Fourier method we have $\|S(t)\|_{H,H}^2\le
e^{-2\lambda_1t}$. We obtain for $n_1<n_2,\,n_1\to\infty$
\begin{align}\label{eqev4}
\begin{split}
\mathbb{E}&\sup_{s\in [0,1]}
\|\eta_1(s+n_1,\theta_{-n_1}\omega)-\eta_1(s+n_2,\theta_{-n_2}\omega)\|_H^2\\
&= \mathbb{E}\sup_{s\in[0,1]}
\|A\int_0^{n_1+s}S(n_1+s-\tau)d\tilde G W(\tau,\theta_{-n_1}\omega)\\
&- A\int_0^{n_2+s}S(n_2+s-\tau)d\tilde G
W(\tau,\theta_{-n_2}\omega)\|_H^2
\end{split}
\end{align}
We can continue by a simple integral substitution, the $(\theta_t)_{t\in\mathbb{R}}$-invariance
of $\mathbb{P}$
and the semigroup property of $S$
\begin{align*}
 &= \mathbb{E}\sup_{s\in [0,1]}
\|A\int_{-n_1}^sS(s-\tau)d\tilde G W(\tau,\omega)-A\int_{-n_2}^sS(s-\tau)d\tilde G W(\tau,\omega)\|_H^2\\
&= \mathbb{E}\sup_{s\in[0,1]}\|A\int_{-n_2}^{-n_1}S(s-\tau)d\tilde
G W\|_H^2=
\mathbb{E}\sup_{s\in[0,1]}\|A\int_{-n_2+n_1}^{0}S(n_1+s-\tau)d\tilde GW\|_H^2\\
&\le \sup_{s\in[0,1]}e^{-2\lambda_1(n_1+s)}\mathbb{E}\|A\int_{-n_1+n_2}^0S(-\tau)d\tilde G W(\tau,\omega)\|_H^2\\
&\le e^{-2\lambda_1n_1} \int_0^\infty\|AS(\tau)\tilde
GQ^\frac{1}{2}\|_{\mathcal{L}_2(U,H)}^2\le Ce^{-2\lambda_1 n_1}.
\end{align*}
The finiteness of the last integral
follows by Lemma \ref{lv1}
Hence $\eta_1(\cdot+n,\theta_{-n}\omega)$ is $L_2$-convergent
with respect to $C([0,1];H)$. This convergence is
exponentially fast.
Via Lemma \ref{lv2} and the semigroup property the exponential convergence to zero of
\begin{equation}\label{last1}
\mathbb{E}\sup_{s\in [0,1]}\|A\int_{-n}^{-n+s}S(s-\tau)d\tilde GW(\tau,\omega)\|_H^2
\end{equation}
for $n\to\infty$ is easily seen.\\
By the Borel and Cantelli Lemma there exists a
set $\Omega^0\subset\Omega$ of full measure such that for
$\omega\in\Omega^0$ we have the above convergences
$\omega$-wise. The same remains true on
$\Omega^j$ if we replace $\omega$ by $\theta_j\omega$ for
$j\in\mathbb{Z}$. The set $\bigcap_j\Omega^j$ is $(\theta_t)_{t\in
\mathbb{Z}}$-invariant and
has a full measure.\\
We denote by $\eta_1(\omega)$ the $\omega$-wise limit of
$\eta_1(0+n,\theta_{-n}\omega)$ which exists on
$\bigcap_j\Omega^j$.
To see that $\eta_1(0+n,\theta_{-n}\theta_{\bar t}\omega)$ is converging on
$\bigcap_j\Omega^j$ we can restrict ourselves to $\bar t\in[0,1]$.
In particular, if $\bar t\in (-1,0)$, $\theta_{\bar
t}\omega$ can be written as
$\theta_{1+\bar t}\theta_{-1}\omega,\;1+\bar t\in (0,1)$ for $\bar t\in (-1,0)$.
We can write
\begin{align*}
\eta_1(n,\theta_{-n+\bar t}\omega)&=A\int_{-n+\bar t}^{\bar t} S(\bar t-\tau)d\tilde GW(\tau,\omega)\\
&=
A\int_{-n}^{\bar t} S(\bar t-\tau)d\tilde GW(\tau,\omega)-
A\int_{-n}^{-n+\bar t} S(\bar t-\tau)d\tilde GW(\tau,\omega).
\end{align*}
By the above considerations the first term on the right hand side is convergent and the second term tends
to zero by (\ref{last1}) using the semigroup property.
Let $\bar\Omega$ be the set of $\omega$'s such that $\eta_1(n,\theta_{-n}\omega)$ converges
$\omega$-wise.
We have for $t\in\mathbb{R}$ that
$\theta_t\bar\Omega\supset\bigcap_j\Omega^j$. Hence the
$(\theta_t)_{t\in\mathbb{R}}$-invariant set
$\bigcap_{t\in\mathbb{R}}\theta_t\bar\Omega$ contains a set of
measure one. We can construct a new probability space if we consider
the trace of $\mathcal{F}$ and the restriction of $\mathbb{P}$ to
this new set $\bigcap_t\theta_t\bar \Omega$. We will use for this new probability
space the old notations.\\
To see that $t\to\eta_1(\theta_t\omega)$ solves (\ref{eqo1}) we
consider $\eta_1(\theta_{t}\omega)-S(t)\eta_1(\omega)$ which is a version of
$A\int_0^tS(t-\tau) d\tilde GW$ such that by the variation of constants formulae
\[
S(t)\eta_1(\omega)+A\int_0^tS(t-\tau)d\tilde GW(\tau,\omega)=\eta_1(\theta_t\omega),
\]
see Da Prato and Zabcyk \cite{DaPZab96} (13.2.13).
\end{proof}
\begin{remark}\label{lastr}
{\rm
Using the technique for the proof of the last lemma we can also show
that there exists a $(\theta_t)_{t\in\mathbb{R}}$-invariant set of full measure such that
$t\to\eta_1(\theta_t\omega)$ has paths in $L_{2,loc}(0,\infty;H^2)$.
Since the distribution of $\eta_1$ is a Gau{\ss} measure it can be derived that
the stationary process
$\eta_1(\theta_\cdot)$ has trajectories in $L_{4,loc}(0,\infty;H^2\cap H)$, see Da Prato and
Zabczyk \cite{DaPZab92} Corollary 2.17.
}
\end{remark}

We now introduce the term {\em tempered random variable}.
   A random variable $\xi$ is called {\em tempered} if this random variable  has a
   subexponential growth:
   \begin{equation*}
   \limsup_{t\to\infty,t\in\mathbb{T}}\frac{\log^+(\xi(\theta_t\omega))}{|t|}=0.
   \end{equation*}
  In the case of an ergodic measure $\mathbb{P}$, it is well known that
there is only one alternative, i.e.,   the $\lim\sup$ in the
   last formula is $+\infty$;  see Arnold \cite{Arn98} Page 165,  Proposition 4.1.3.

   \begin{lemma}\label{llll}
   There exists a
   $(\theta_t)_{t\in\mathbb{R}}$-invariant
   set $\Omega_1$ of full $\mathbb{P}$-measure such that for any
   $\omega\in\Omega_1$ , $\eta_1(\omega)$ is well defined and
   \[
   \lim_{t\to\pm\infty}\frac{\log^+\|\eta_1(\theta_t\omega)\|_{H}}{t}=0
   \]
   which means that the $H$ norm of $\eta_1$ is tempered.
   \end{lemma}
   \begin{proof}
   On account of (\ref{eqev4}) we have
   \[
   \mathbb{E}\sup_{s\in
   [0,1]}\|\eta_1(\theta_s\omega)\|_{H}<\infty.
   \]
   Hence by the Birkhoff ergodic theorem, we have a
   $(\theta_t)_{t\in\mathbb{Z}}$--invariant set $\Omega_1\subset \Omega$ such that
   \[
   \lim_{n\to\pm\infty}
   \frac{\sup_{s\in[0,1]}\|\eta_1(\theta_{s+n}\omega)\|_{H}}{n}=0.
   \]
   On the other hand for $t>0$:
   \[
   \|\eta_1(\theta_{t}\omega)\|_{H}
   \le
   \sup_{s\in[0,1]}\|\eta_1(\theta_{s+[t]}\omega)\|_{H}
   \]
   such that
   \[
   \lim_{t\to\pm\infty}
   \frac{\|\eta_1(\theta_{t}\omega)\|_{H}}{t}=0.
   \]
   and similarly for $t<0$, which is sufficient for the
   $(\theta_t)_{t\in\mathbb{R}}$--invariance of $\Omega_1$.
   \end{proof}

Collecting our
calculations  we get
\begin{theorem}
Let $Q$ be a linear operator on
$U=H^\frac{5}{2}(O_{\cdot,\cdot,2\pi})\cap
\dot{L}_2(O_{\cdot,\cdot,2\pi})$ with finite trace. Then there
exists a tempered random variable $\eta_1$ and a
$(\theta_t)_{t\in\mathbb{R}}$-invariant set of full measure such
that
\[
t\to\eta_1(\theta_t\omega)\in C([0,\infty);H)\cap L_{4,loc}(0,\infty;H^2\cap H).
\]
$\eta_1(\omega)$ is periodic in $x,\,y$ directions.
\end{theorem}

   In a similar manner we can consider (\ref{eqo1}) with time-periodic
   boundary condition representing the impact of the earth's rotation on the fluid :
   \begin{equation}\label{eqo3}
   {\eta_2}_t=\tilde\Delta \eta_2,\;\frac{\partial\eta_2}{\partial
   n}=0,\;\text{on }
   O_{\cdot,\cdot,0},
   \quad \frac{\partial\eta_2}{\partial n}=u_0\sin(2\pi
   t)\;\text{on }
   O_{\cdot,\cdot,2\pi}
   \end{equation}
   for $u_0\in H^{\frac{1}{2}}\cap \dot{L}_2(O_{\cdot,\cdot,2\pi})$ and (spatial) periodic boundary
   conditions on the other faces. We obtain without proof the following :
   \begin{lemma}
   Suppose that $u_0\in H^\frac{1}{2}(O_{\cdot,\cdot,2\pi})\cap\dot{L}_2(O_{\cdot,\cdot,2\pi})$. Then there exists a
   continuous periodic
   solution
   \[
   t\to\eta_2(t)\in H^2
   \]
   which satisfies (\ref{eqo3}). In particular, $\eta_2(t)$ is
   also periodic in $x$-$y$ directions.
   \end{lemma}
   The proof is based on the properties of the operator $\tilde G$.\\

   Let $\theta^2=(\theta^2)_{t\in\mathbb{R}}$ be the shift operator
   $\theta_t^2f(\cdot)=f(t+\cdot)$ for appropriate functions $f$.
   We consider the hull of $u_0\sin t$ with respect to $\theta^2$:
   \[
   \Omega_2=\bigcup_{t\in\mathbb{R}}\theta_t^2 (u_0\sin(2\pi\cdot))=
   \bigcup_{t\in\mathbb{R}}u_0\sin(2\pi(t+\cdot))
   =\bigcup_{t\in[0,2\pi)}u_0\sin(2\pi(t+\cdot))
   .
   \]

   Summarizing, we have found a process $\eta=\eta_1+\eta_2$ which will serve
   as a model for the perturbation on the ocean surface.\\

    After these preparations we can introduce a nonautonomous/random
   dynamical system. Let $\theta$ be a flow on a set
   $\Omega$ (such that (\ref{eqo6}) is
   fulfilled).
   \begin{definition}
   A nonautonomous dynamical system $\phi$ on a phase space $H$ with respect
   to a flow $\theta$
   is a mapping
   \[
   \phi:\,
   \mathbb{T}^+\times\Omega\times H\to H
   \]
   fulfilling the cocycle property
   \begin{align*}
   \phi(t+\tau,\omega,\cdot)&=\phi(t,\theta_\tau\omega,\phi(\tau,\omega,\cdot))
   \quad\text{for }t,\,\tau\in\mathbb{T}^+\\
   \phi(0,\omega,\cdot)&={\rm id}_H
   \end{align*}
   for $\omega\in\Omega$ and $t,\,\tau\in\mathbb{T}^+$. Suppose that
   the flow $\theta$ is carried by a metric dynamical system and
   $\phi$ is supposed to be measurable then $\phi$ forms a random dynamical
   system.
   \end{definition}
   This definition is due to Arnold \cite{Arn98} or Kloeden et. al. \cite{KloKelSchm97a}
   \begin{remark}\label{rdis}
   {\rm
   We now consider the flow $\theta=(\theta^1,\theta^2)$ on
   $\Omega=\Omega^1\times\Omega^2$.
   By the global forward existence and uniqueness of the solution of
   (\ref{eq2.2}) for
   $\eta(\theta_t\omega)=\eta_1(\theta^1_t\omega_1)+\eta_2(\theta^2_t\omega_2),\,
   \omega=(\omega_1,\omega_2)\in\Omega$, the solution operator of (\ref{eq2.2}), which maps
   an initial condition $u_0\in H$ and a  sample point  $\omega$
   to the solution at time $t$, has the cocycle property for $H=\dot{L}_2(O)$.
   We will denote this
solution operator
     by $\phi(t,\omega,x)$. Note that by the periodicity of $\eta_2$
   the restriction of $\phi$ on $\mathbb{Z}$ is a {\em random dynamical
   system} for any fixed $\omega_2\in\Omega_2$.
   Indeed, $\theta_i\omega=(\theta_i^1\omega_1,\omega_2)$ which can be identified with
   $\theta^1$  and which leaves $\mathbb{P}$
   invariant ($\mathbb{P}$ is ergodic with respect to
   $(\theta_i^1)_{i\in\mathbb{Z}}$).
   }
   \end{remark}
   \begin{remark}
   {\rm
   Another opportunity to get an example for a {\em complete} random
   dynamical system would be to equip $\Omega_2$ with an ergodic
   measure. But in contrast to the fact that the daily or yearly
   rotation of the earth is well determined such a random ansatz
   would express that the beginning of these periods is rather
   random.
   }
   \end{remark}
   The main result of this article is   the existence of an attractor
   for the nonautonomous dynamical system. This attractor will attract
   random sets in probability.  Before we give the main theorem, we
   make
   some basic remarks on {\em random sets}. \\
   Suppose that $H$ is a Polish space. A set function $\omega\to D(\omega)$
   with closed and nonempty images is called a closed random set over
   $(\Omega_1,\mathcal{F},\mathbb{P})$
   if and only if there exists a countable number of random variables
   \[
   \xi_i:(\Omega_1,\mathcal{F},\mathbb{P})\to H,\quad i\in\mathbb{N}
   \]
   such that
   \[
   D(\omega)=\overline{\bigcup_{i\in\mathbb{N}}\xi_i(\omega)}
   \]
   see Castaing and Valadier \cite{CasVal77} Chapter 3.

   A random set $D$ is called tempered if the random variable
   ${\rm
   dist}_H(D(\omega),\{0\})$ is tempered
   where
   \[
   {\rm dist}_H(A,B)=\sup_{a\in A}\inf_{b\in B}\|a-b\|_H.
   \]

   We now define the term {\em random attractor}.
   \begin{definition}
   Let $\phi$ be a random dynamical system over the metric
   dynamical system $\theta$. A tempered random set
   $(A(\omega))_{\omega\in\Omega}$ with compact and nonempty
   images $A(\omega)$ is called random attractor if
   \begin{equation}\label{eq20}
   \phi(t,\omega,A(\omega))=A(\theta_t\omega),\qquad t\in
   \mathbb{T}^+
   \end{equation}
   and for any random tempered closed random set $D(\omega)$:
   \begin{equation}\label{eq21}
   (\mathbb{P})\lim_{t\to\infty,t\in\mathbb{T}^+}{\rm dist}_H
   (\overline{\phi(t,\omega,D(\omega))},A(\theta_t\omega))=0.
   \end{equation}
   \end{definition}
   In the next section we will show that the dynamical system
   generated by (\ref{eq2.2}) restricted to $\mathbb{T}=\mathbb{Z}$
   has such a random attractor. However, because of   time-periodic
   perturbation,   we do not have a random dynamical system but a nonautonomous
   dynamical system for  $\mathbb{T}=\mathbb{R}$.
   Therefore we have to modify this conclusion for
   $\mathbb{T}=\mathbb{Z}$ a little bit, and thus we obtain
the following main result in this section.

   \begin{theorem} \label{maintheorem}
   There exists set function  $A=(A(\omega))_{\omega\in\Omega}$
for the 3D QGE  (\ref{eq2.2}) under random
plus
 time-periodic forcing on the ocean surface such that (\ref{eq20})
(\ref{eq21})is satisfied for set functions $D$ with closed images
and subsexponential growth of $t\to {\rm
dist}_H(D(\theta_t\omega),\{0\})$ with tempered restriction to
$(\theta_t)_{t\in\mathbb{Z}}$. Moreover,
   the mapping
   \[
   t\to A(\theta_t\omega)
   \]
   has a subexponential growth.
   \end{theorem}

In the next section, we prove this main theorem \ref{maintheorem}.


   \section{Proof of the main theorem}

   The proof of the main theorem \ref{maintheorem} in the last section
is based on checking following conditions.
   We will consider the set  of families of sets $D$ such
   that $D(\omega)$ is closed and nonempty and $t\to {\rm dist}_H(D(\theta_t(\omega_1,\omega_2)),0)$
   has a subexponintially growth. If we consider $D$ with respect to the
   restriction of  $(\theta_t)_{t\in\mathbb{R}}$ to $\mathbb{Z}$ family of sets $D$ is supposed to
   be tempered. The set of these set families  is denoted by
   $\mathcal{D}$.\\

   The following theorem states the existence of random attractors.
   \begin{theorem}\label{t2}
   Suppose that $\theta^1$ be a metric dynamical system and
   suppose that for $\omega_1\in\Omega_1,\,t\in\mathbb{Z}^+$ the mappings
   $\phi(t,\omega_1,\cdot)$ are continuous and that there exists
   a tempered set  $B$ having compact images which
   is absorbing:
   \begin{equation}\label{eqx100}
   \phi(t,\theta_{-t}^1\omega_1,D(\theta_{-t}^1\omega_1))\subset B(\omega_1)
   \end{equation}
   for $\mathbb{Z}^+\ni t\ge t_0(D,\omega)$, $D$ is supposed to be tempered.
   Then there exists a unique random attractor $A$ which is tempered.
   \end{theorem}

   The proof of this theorem can be found in Flandoli and Schmalfu{\ss}
   \cite{FlaSchm95a}.\\

   Now let us use this result to prove our main Theorem
   \ref{maintheorem}.

\begin{proof}

   Calculating the inner product  in $H$ we obtain by Lemma
   \ref{c1}, \ref{ll}, \ref{ll1}.
   \begin{equation}
   \begin{split}
   \|u(t)\|_H^2&+2\nu \int_0^t\|u(\tau)\|_V^2d\tau=\|u_0\|_H^2
   +2\beta\int_0^t\langle G(\tilde \Delta\eta)_x-\eta_x,u\rangle d\tau\\
   &\le\|u_0\|_H^2+\frac{\beta^2}{\nu}\int_0^t\|G(\tilde\Delta\eta(\theta_\tau\omega))_x-\eta(\theta_\tau\omega)_x\|_{V^\prime}^2d\tau+
   \nu\int_0^t\|u(\tau)\|_V^2 d\tau.\label{eqx2}
   \end{split}
   \end{equation}
   Parallel to this inequality we consider the affine equation
   \begin{equation}\label{x1}
   \begin{split}
   \xi(t)&+\nu \lambda_1
   \int_0^t\xi(\tau)d\tau\\
   &
   =\|u_0\|_H^2+
   \frac{\beta^2}{\nu}\int_0^t\|G(\tilde\Delta\eta(\theta_\tau\omega))_x-\eta(\theta_\tau\omega)_x\|_{V^\prime}^2d\tau.
   \end{split}
   \end{equation}
   The solution of this equation $\xi(t,\omega,\|u_0\|_H^2)$ is a bound for
   $\|u(t)\|^2_H$. In addition,
   this equation has a unique forward and backward exponentially fast
   attracting solution
   \[
   (\omega,t)\to\xi^\ast(\theta_t\omega),\quad
   \xi^\ast(\omega)=\frac{\beta^2}{\nu}\int_{-\infty}^0
   e^{\nu \lambda_1\tau}
   \|G(\tilde\Delta\eta(\theta_\tau\omega))_x-\eta(\theta_\tau\omega)_x\|_{V^\prime}^2d\tau.
   \]
   for $t\to\infty$ which follows by the variation of constants formulae.\\
   To see the backward convergence uniform with respect to $D$
   we must prove that
   \[
   \lim_{t\to\infty}\sup_{u_0\in
   D(\theta_{-t}\omega)}|\xi(t,\theta_{-t}\omega,\|u_0\|_H^2)-\xi^\ast(\omega)|=0.
   \]
   Because of
   \[
   \xi^\ast(\omega)=\xi(t,\theta_{-t}\omega,\xi^\ast(\theta_{-t}\omega)))
   \]
   we have
   \[
   \sup_{u_0\in
   D(\theta_{-t}\omega)}|\xi(t,\theta_{-t}\omega,\|u_0\|_H^2)-\xi^\ast(\omega)|
   \le
   \sup_{u_0\in
   D(\theta_{-t}\omega)}|\|u_0\|_H^2-\xi^\ast(\theta_{-t}\omega)|e^{-\lambda_1t}
   \]
   which also follows by the variation of constants formulae.
   Since the first factor on the right hand side is only subexponentially growing the
   convergence conclusion follows for $t\to\infty$.\\
   The  forward convergence follows
   straightforwardly by the variation of constants formulae
   \[
   |\sup_{u_0\in
   D(\omega)}\xi(t,\omega,\|u_0\|_H^2)-\xi^\ast(\theta_t\omega)|\le
   e^{-\lambda_1t}\sup_{u_0\in D(\omega)}|\|u_0\|_H^2-\xi^\ast(\omega)|
   \]
   where $\xi^\ast(\theta_t\omega)=\xi(t,\omega,\xi^\ast(\omega))$.
   The mapping $t\to\xi^\ast(\theta_t\omega)$ is subexponentially
   growing. Indeed, we have
   \begin{align*}
   \frac{\nu}{2\beta^2}\xi^\ast(\omega)&\le
   \int_{-\infty}^0
   e^{\nu \lambda_1\tau}
   \|G(\tilde\Delta\eta_1(\theta_\tau\omega))_x-\eta_1(\theta_\tau\omega)_x\|_{V^\prime}^2d\tau\\
   &+
   \int_{-\infty}^0
   e^{\nu \lambda_1\tau}
   \|G(\tilde\Delta\eta_2(\theta_\tau\omega))_x-\eta_2(\theta_\tau\omega)_x\|_{V^\prime}^2d\tau.
   \end{align*}
   The second term on the right hand side is bounded and
   the first is tempered, see Theorem 4.1 in \cite{Schm97a}. Since $\|\eta_1\|_H$
   is tempered
   so is $\|G(\tilde\Delta\eta_1)_x\|_{V^\prime}$\\
   Hence the ball $\tilde B$ with center zero and radius $2\xi^\ast(\omega)$ forms
   an absorbing set (in the sense of (\ref{eqx100})), $\tilde
   B(\omega)\in\mathcal{D}$.
   Plugging in the radius  $2\xi^\ast(\omega)$ of $\tilde B$ into (\ref{x1}) for $\|u_0\|_H^2$ it is
   easily
   seen that $\tilde B$ is forward invariant:
   \[
   \phi(t,\omega,\tilde B(\omega))\subset \tilde B(\theta_t\omega).
   \]
   Since for fixed $\omega_2$ the mapping $\xi^\ast$ is a random variable $B$ is a random set.
   On the other hand we can check that
   \[
   B(\omega):=\overline{\phi(1,\theta_{-1}\omega,\tilde
   B(\theta_{-1}\omega))}\subset \tilde B(\omega).
   \]
   Note that $B$ is compact because $\phi$ is regularizing
   which follows by Remark \ref{rneu} ii).
   We have seen that $B$ is forward and backward absorbing.
   This remains true if we restrict our random dynamical system
   to discrete time, see Remark \ref{rdis}. By the above inclusion we can
   derive $B\in \mathcal{D}$ since $\tilde B\in\mathcal{D}$. It follows also by the definition
   of $B$ that for fixed $\bar \omega_2$ the set $B((\omega_1,\bar\omega_2))$ is random
   set which is tempered.
   \\

   The continuity of $\phi(t,\omega,\cdot)$ follows
   because we have the estimate from Lemma \ref{l3}. However, we need this
   technique
   later on once more such that we are going to explain this technique:
   We have by Lemma \ref{l3} and the properties of $C(t)$ and $D$:
   \begin{eqnarray*}
   &&
   |\langle B(u_1,u_1)-B(u_2,u_2),u_1-u_2\rangle|\le c_3
   \|u_1-u_2\|_H\|u_1-u_2\|_V\|u_1\|_V\\
   &&
   \langle C(t,u_1)-C(t,u_2),u_1-u_2\rangle=0\\
   &&
   \langle D(u_1)-D(u_2),u_1-u_2\rangle=0
   \end{eqnarray*}
   Hence there exists a constant $c_{6}$ such that
   for two solutions $u_1,\,u_2$ with initial conditions $u_{10},\, u_{20}$
   \begin{equation*}
   \frac{d}{dt}\|u_1(t)-u_2(t)\|_H^2
   \le c_{6}(\|u_1(t)\|_V^2)
   \|u_1(t)-u_2(t)\|_H^2.
   \end{equation*}
   It follows by the Gronwall lemma
   \begin{equation}\label{eqx101}
   \|u_1(t)-u_2(t)\|_H^2\le \|u_{10}-u_{20}\|_H^2
   e^{c_6\int_0^t
   \|u_1(\tau)\|_V^2d\tau}
   \end{equation}
   which gives the continuity of $\phi(t,\omega,\cdot)$.
    \\
   We have checked all assumptions of the above theorem such that
   the discrete random dynamical system has a random attractor
   $A(\omega_1,\bar\omega_2)$ for some fixed $\bar\omega_2$ with respect to the restricted
   discrete random dynamical system.
   We now extend the definition of $A$ to $\Omega_2$. We set:
   \[
   A(\theta_t(\omega_1,\bar\omega_2))=\phi(t,(\omega_1,\bar\omega_2),A(\omega_1,\bar\omega_2))\qquad
   \text{for }t\in\mathbb{R}^+
   \]
   which is invariant in the sense of (\ref{eq20}).
   Note that by the cocycle property this definition is correct.
   We show convergence (\ref{eq21}):
   \[
   (\mathbb{P})\lim_{t\to\infty,t\in\mathbb{R}^+}
   {\rm dist}_H(\overline{\phi(t,(\omega_1,\bar\omega_2),
   D((\omega_1,\bar\omega_2)))},A(\theta_t(\omega_1,\bar\omega_2)))=0.
   \]
   for $D\in\mathcal{D}$ and any $\bar\omega_2\in\Omega_2$.
   Since $B(\omega)$ is a forward absorbing
   and forward invariant set for any $D\in\mathcal{D}$
   it remains to check the convergence conclusion for $D= B$.
   On account of (\ref{eqx2})
   we can notice for fixed $\bar\omega_2$ and
   $u_1=\phi(t,(\omega_1,\bar\omega_2),y)$
   \begin{align*}
   \sup_{y\in B((\omega_1,\bar\omega_2))}&\nu
   \int_0^1\|\phi(\tau,(\omega_1,\bar\omega_2),y)\|_V^2d\tau\\
   &
   \le 2\xi^\ast(\omega_1,\bar\omega_2)+
   \frac{\beta^2}{\nu}\int_0^1
   \|G(\tilde\Delta\eta(\theta_\tau(\omega_1,\bar\omega_2)))_x
   -\eta(\theta_\tau(\omega_1,\bar\omega_2))_x\|_{V^\prime}^2d\tau.
   \end{align*}
   With respect to (\ref{eqx101}) we have  an appropriate constant $c_{7}>0$:
   \begin{align*}
   &\sup_{x\in B((\omega_1,\bar\omega_2))}e^{c_6\int_0^t\|\phi(\tau,(\omega_1,\bar\omega_2),x)\|_V^2d\tau}
   \\
   &\quad\le
   e^{c_{7}(\xi^\ast(\omega_1,\bar\omega_2)+\int_0^1
   \|G(\tilde\Delta\eta(\theta_\tau(\omega_1,\bar\omega_2))_x
    -\eta(\theta_\tau(\omega_1,\bar\omega_2))_x\|_{V^\prime}^2 d\tau)}
   =:Y((\omega_1, \bar\omega_2))
   \end{align*}
   for $t\in[0,1]$.
   The mapping
   \[
   (\omega_1,n)\to Y((\theta_{n}\omega_1,\bar\omega_2))
   \]
   defines a stationary process. We obtain by (\ref{eqx101}):
   \begin{align*}
   \sup_{\tau\in[0,1]}&{\rm
   dist}_H(\phi(\tau+n,\omega,B((\omega_1,\bar\omega_2)),
   A(\theta_{\tau+n}(\omega_1,\bar\omega_2)))\\
   &
   \le
   {\rm dist}_H(\phi(n,(\omega_1,\bar\omega_2),B((\omega_1,\bar\omega_2)),
   A((\theta_n^1\omega_1,\bar\omega_2)))Y((\theta_n^1\omega_1,\bar\omega_2))^\frac{1}{2}.
   \end{align*}
   Since the first factor of the right hand side tends to zero in probability
   the product of the right hand side also tends to zero in probability
   which gives the convergence conclusion (\ref{eq21}).\\
   To obtain the subexponential growth of
   \[
   t\to A(\theta_t(\omega_1,\bar\omega_2))
   \]
   we need that $Y$ is tempered with respect to the restricted flow
   $(\theta_t)_{t\in\mathbb{Z}}$.
   Due to Arnold \cite{Arn98} Proposition 4.1.3 it is sufficient to
   show that
   \[
   \mathbb{E}\xi^\ast(\omega_1,\bar\omega_2)<\infty,\quad\mathbb{E}
    \int_0^1\|G(\tilde\Delta\eta(\theta_{\tau+s}(\omega_1,\bar\omega_2)))_x
    -\eta(\theta_{\tau+s}(\omega_1,\bar\omega_2))_x\|_{V^\prime}^2d\tau<\infty
   \]
   which can be derived from Lemma \ref{lv2}.
   Consequently,
   \[
   \sup_{x\in A(\theta_t(\omega_1,\bar\omega_2))}\|x\|_H\le Y((\omega_1,\bar\omega_2))^\frac{1}{2}{\rm
   diam}(A((\omega_1,\bar\omega_2)))+\sup_{x\in A((\omega_1,\bar\omega_2))}\|x\|_H,\quad t\in[0,1]
   \]
   by the triangle inequality. We have used that the
   product of two tempered random variables is tempered
   which gives the general convergence conclusion of the  main theorem \ref{maintheorem}
   since $A(\omega_1,\bar\omega_2)$ is $(\theta_t)_{t\in \mathbb{Z}}$-tempered.
\end{proof}


\section{Summary}

We have studied the 3D  baroclinic quasigeostrophic flow model   under
random wind forcing and time-periodic fluctuations on fluid
boundary, i.e.,  on the interface between the oceans and the atmosphere.
The time-periodic fluctuations are due to periodic rotation
of the earth and thus periodic exposure of the earth to the solar
radiation.
We have established the  well-posedness of the baroclinic
quasigeostrophic flow model in the state space (Theorem \ref{tex}),
and we have demonstrated the
existence of the random attractors (Theorem \ref{maintheorem}), again in the state space.  We
   have also discussed the relevance of our results
 to climate modeling.

\bigskip

{\bf Acknowledgement.} A part of this work was done at the
Oberwolfach Mathematical Research Institute supported by {\em
Volkswagen Stiftung} while the authors were Research in Pairs
Fellows. This work was partly supported by the NSF Grant
DMS-9973204.

    \end{document}